\documentclass[12pt]{amsart}
\usepackage{amsfonts}

\usepackage{amscd}

\newtheorem{proposition}{Proposition}[section]
\newtheorem{theorem}[proposition]{Theorem}
\newtheorem{lemma}[proposition]{Lemma}
\newtheorem{corollary}[proposition]{Corollary}

\newtheorem{definition}[proposition]{Definition}

\newtheorem{example}[proposition]{Example}
\begin{document}
\author{Emmanuel Breuillard}
\date{\today}
\title{Equidistribution of dense subgroups on nilpotent Lie groups}
\address{Emmanuel Breuillard, Ecole Polytechnique, 91128 Palaiseau, France}
\email{emmanuel.breuillard@math.polytechnique.fr}

\begin{abstract}
Let $\Gamma $ be a dense subgroup of a simply connected nilpotent Lie group $%
G$ generated by a finite symmetric set $S.$ We consider the $n$-ball $S^{n}$
for the word metric induced by $S$ on $\Gamma $. We show that $S^{n}$ (with
uniform measure) becomes equidistributed on $G$ with respect to the Haar
measure as $n$ tends to infinity. We also prove the analogous result for
random walk averages.
\end{abstract}
\maketitle

\section{Introduction}

In this paper we are concerned with equidistribution properties of dense
subgroups of Lie groups. Let $G$ be a Lie group and $S=%
\{1,a_{1},...,a_{s},a_{1}^{-1},...,a_{s}^{-1}\}$ a finite symmetric set of
elements in $G$, which we assume to generate a dense subgroup $\Gamma $ in $G
$. Let $S^{n}$ be the $n$-ball in the Cayley graph of $\Gamma $ induced by $S
$, i.e. the set of elements of $\Gamma $ that can be written as a product of
at most $n$ elements from $S$. We give ourselves two open subsets $U$ and $V$
in $G$ with Lebesgue-negligible boundary$.$ Consider the ratios 
\begin{equation*}
R_{n}(U,V)=\frac{|S^{n}\cap U|}{|S^{n}\cap V|}
\end{equation*}
where $|A|$ denotes the cardinal of a set $A$. The question is to find out
whether or not $R_{n}(U,V)$ converges and if it does to determine the
measure $m$ on $G$ such that 
\begin{equation}
\lim_{n\rightarrow +\infty }R_{n}(U,V)=\frac{m(U)}{m(V)}  \label{ratio}
\end{equation}

Various kinds of averages are possible instead of the uniform averages over $%
S^{n}$ above. In particular one may consider random walk averages or
averages over the abstract free group generated by $S$. Limits such as (\ref
{ratio}) and results of this kind are usually called \textit{ratio limit
theorems} (see \cite{LeP}, \cite{Guiv2}).

Arnol'd and Krylov were among the first to consider such a question and they
showed in \cite{ArK} how to use the representation theory of $G$ to quickly
give a positive answer to this problem when $G$ is a compact Lie group and $%
S $ is the set of free generators of a dense free subgroup of rank $s$ in $%
G. $ Then (\ref{ratio}) holds with $m$ the normalized Haar measure on $G$.

In this paper we prove that (\ref{ratio}) holds with $m$ a Haar measure when 
$G$ is any closed subgroup of a simply connected nilpotent Lie group
(Corollary \ref{coro} below). The case of semi-simple Lie groups (even $%
SL_{2}(\Bbb{R})$) remains an open question up to now. In \cite{Kaz}, Kazhdan
considered the case when $G=Isom(\Bbb{R}^{d})^{\circ }$ for $d=2$ and
obtained an analogous ratio limit theorem for random walk averages. His
argument was corrected and the result extented by Guivarc'h in \cite{Guiv2}
(see also \cite{Bre} for the precise asymptotics) but the case when $d>2$ is
a well-known open problem.

Let us fix some notation and then state the main result of this paper. For a
group $G$, we let $vol_{G}$ be a left Haar measure on $G$. Recall that
according to a theorem of Guivarc'h (see \cite{Guiv} and \cite{Bre2}), if $G$
is nilpotent, locally compact and compactly generated, then there is an
integer $d(G)\in \Bbb{N}$ such that $vol_{G}(\Omega ^{n})\approx n^{d(G)}$
for any compact generating set $\Omega $. When $G$ is connected and simply
connected, the integer $d(G)$ is given by the Bass-Guivarc'h formula (\ref
{exp}) below. We show:

\begin{theorem}
\label{main}(Dense subgroups are equidistributed) Let $\Gamma $ be a
finitely generated nilpotent group and $S=%
\{1,a_{1},...,a_{s},a_{1}^{-1},...,a_{s}^{-1}\}$ a finite generating set.
Let $G$ a closed subgroup of a simply connected nilpotent Lie group. Suppose 
$\phi :\Gamma \rightarrow G$ is a homomorphism with dense image. Then there
is a positive constant $C>0$ such that 
\begin{equation}
\lim_{n\rightarrow +\infty }\frac{|S^{n}\cap \phi ^{-1}(B)|}{n^{d(\Gamma
)-d(G)}}=C\cdot vol_{G}(B)  \label{a}
\end{equation}
for every bounded Borel subset $B\subset G$ with negligible boundary.
\end{theorem}

\begin{corollary}
\label{coro}If $G$ is a closed subgroup of a simply connected Lie group,
then (\ref{ratio}) holds for $m$ a Haar measure on $G$.
\end{corollary}

This result can be seen as a generalization of the classical Weyl
equidistribution (\cite{Wey}) of multiples of an irrational number $\alpha $
modulo $1$. Indeed, let $G=\Bbb{R}$, $\Gamma =\Bbb{Z}^{2}$, and $%
S=\{(0,0),(0,\pm 1),(\pm 1,0)\}$, $\phi (x,y)=x+\alpha y$ and $%
B=(a,b)\subset [0,1].$ Then (\ref{a}) in this case translates as 
\begin{equation*}
\lim_{n\rightarrow +\infty }\frac{1}{2n}|\{k\in [-n,n],k\alpha \in (a,b)%
\text{ mod }1\text{ }\}|=b-a.
\end{equation*}

For abelian $G$ the theorem can be deduced from Weyl's equidistribution.
However, as we will see below, when $G$ is nilpotent, it requires different
techniques.

Combining Alexopoulos' theorem \cite{Alex} on the asymptotics of the return
probability of random walks on finitely generated nilpotent groups with
Theorem \ref{main} allows to prove the analogous limit theorem for random
walk averages, namely:

\begin{corollary}
\label{coro2}(Local Limit Theorem) Let $G$ be a simply connected nilpotent
Lie group and $\mu $ a symmetric and finitely supported and probability
measure on $G$ whose support generates a dense subgroup, then there is $%
c(\mu )>0$ such that for any Borel set $B$ with Lebesgue-negligible
boundary, 
\begin{equation}
\lim_{n\rightarrow +\infty }n^{-d(G)/2}\mu ^{*n}(B)=c(\mu )\cdot vol_{G}(B)
\label{b}
\end{equation}
where $\mu ^{*n}$ denotes the $n$-th fold convolution power of $\mu $.
\end{corollary}

The proof of Theorem \ref{main} makes use of two crucial ingredients. First,
we need precise information on the shape of the $n$-balls $S^{n}$ in $\Gamma
,$ and this is essentially provided by Pansu's theorem from \cite{Pan}.
Second, we use a now well-known principle from ergodic theory (see \cite{DRS}%
, \cite{Led}) according to which the ergodic properties of the action of $%
\Gamma $ on a homogeneous space $N/M$ are dual to those of the action of $M$
on $N/\Gamma $. The unique ergodicity of unipotent flows on nilmanifolds,
which is an ``ancestor'' of Ratner's theorem, allows then to reduce the
equidistribution statement to a computation of the asymptotic volume of
cosets of $M$ inside large balls in $N$. The goal of Section 2 is to prove
this volume estimate (Proposition \ref{ThmDiscs}) and give some background
on homogeneous quasi-norms on nilpotent Lie groups. In Section 3, we
complete the proof of Theorem \ref{main}.

\smallskip

\textbf{N.B.:} \textbf{a. }We prove (see Theorem \ref{eq} and \ref{unifdist}) that the limit (%
\ref{a}) exists also for more general averages of the form $B(n)\cap \phi
^{-1}(U)$ in place of $S^{n}\cap \phi ^{-1}(U)$, where $B(n)$ is the $n$%
-ball for a quasi-norm on $\Gamma $, or any left-invariant coarsely geodesic
distance.

\textbf{b.} The techniques of this paper allow to get uniformity of
convergence in (\ref{a}) (resp. (\ref{b})), when $B$ is allowed to vary
among translates $xB$ such that the distance between $x$ and $1$ is a $o(n)$
(resp. $o(\sqrt{n})$)$,$ and a uniform upper bound for the left hand side of
(\ref{a}) holds for all translates. But we will not need these refinements
here.

\textbf{c.} Our equidistribution problem of a dense subgroup $\Gamma $ in a
nilpotent Lie group $G$ can be phrased more generally as the question of
whether $\Gamma $-orbits equidistribute in their closure on a homogeneous
space $N/M.$ Here we treated the case of $M$ normal in $N$ and a dense
orbit. The general case, when $N$ is nilpotent, is slightly more involved
but can be treated by similar methods.

\textbf{d.} In \cite{Bre3} we showed, in the case when $G$ is the Heisenberg
group, that (\ref{b}) holds for more general measures, namely any centered
and compactly supported measure $\mu $ on $G.$ The methods of \cite{Bre3}
use representation theory and are of a very different nature from the proof
displayed in the present paper.

\tableofcontents%

\section{Homogeneous quasi-norms and volume of balls on nilpotent Lie groups}

The main goal of this section is to prove the following result for a simply
connected nilpotent Lie group $N.$

\begin{proposition}
\label{propdeb}For any quasi-norm on $N,$ the balls $(D_{t})_{t>0}$, $%
D_{t}=\{|x|\leq t\},$ form a nicely growing family of subsets of $N$.
\end{proposition}

Quasi-norms on $N$ and nicely growing subsets are defined below. This
statement essentially means that given a left invariant metric on $N$, we
can compute the asymptotics of the volume of balls for the induced metric on
a closed connected subgroup. To prove it, we will need to introduce some
background on nilpotent Lie groups. This is also the second goal of this
section. In the next two paragraphs we deal exclusively with filtered vector
spaces, while in the remaining four we apply the results to nilpotent Lie
groups.

\subsection{Filtrations on vector spaces and exterior powers}

On a nilpotent Lie algebra, there is a canonical filtration given by the
descending central series. In this paragraph, we describe some properties of
filtrations and associated degree functions in the more general context of
vector spaces. This allows to attach a degree $\deg _{V}(W)$ to every vector
subspace $W$, a notion that will be useful when proving Theorem \ref
{ThmDiscs}. The content of this paragraph is probably well known but we
couldn't find an adequate reference for it.

\subsubsection{Filtrations and degree on $V$ and $\Lambda ^{*}V$}

Let $V$ be a real vector space equipped with a filtration, i.e. a
non-increasing finite sequence of vector subspaces $V=V_{1}\supseteq
V_{2}\supseteq ...\supseteq V_{r}$.

\begin{proposition}
\label{degprop}Associated to this filtration is a function $\deg
:V\rightarrow \Bbb{N}$ called the \textit{degree} defined for $v\in V$ by $%
\deg (v)=\max_{i\geq 1}\{i,v\in V_{i}\}.$ This degree function extends in a
canonical way to the exterior power $\Lambda ^{*}V.$
\end{proposition}

To see this, one can for instance consider a basis $(e_{1},...,e_{n})$ of $V$
which is adapted to the filtration $(V_{i})_{i}$ in the sense that $%
V_{i}=span\{e_{k}|k=1,...,n,$ $\deg (e_{k})\geq i\}.$ This basis gives rise
to an associated basis for $\Lambda ^{*}V$ given by the $e_{I}=e_{i_{1}}%
\wedge ...\wedge e_{i_{k}}$'s where $I=\{i_{1},...,i_{k}\}$ ranges over all
subsets of $\{1,...,n\}$. We can then define the degree of a basis element
by 
\begin{equation}
\deg (e_{I})=\deg (e_{i_{1}})+...+\deg (e_{i_{k}}).  \label{degext}
\end{equation}
Subsequently, this defines a filtration on $\Lambda ^{*}V$ by letting $%
\widehat{V_{i}}=span\{e_{I}|I\subseteq \{1,...,n\},\deg (e_{I})\geq i\}.$ In
turn, we get a degree function on $\Lambda ^{*}V$ extending the definition (%
\ref{degext}) by setting $\deg (\xi )=\max_{i\geq 1}\{i,\xi \in \widehat{%
V_{i}}\}.$ It is just a matter of simple verification to check that the
filtration $(\widehat{V_{i}})_{i}$ and the degree thus defined on $\Lambda
^{*}V$ are independent of the choice of an adapted basis.

\subsubsection{Induced filtration on a subspace, degree of a subspace\label%
{degreesub}}

Let $W$ be a vector subspace of $V.$

\begin{proposition}
The operations of restricting to a subspace $W$ and extending to the
exterior power commute and give rise, after composition, to a uniquely
defined degree function $\deg :\Lambda ^{*}W\rightarrow \Bbb{N}.$
\end{proposition}

\proof%
Given an adapted basis of $W$ with respect to the filtration $(W\cap
V_{i})_{i},$ it is possible to complete this basis into a basis $%
(e_{1},...,e_{n})$ of $V$ which is adapted to the original filtration $%
(V_{i})_{i}$. So we easily check that 
\begin{equation*}
\widehat{W\cap V_{i}}=span\{e_{I}|I\subseteq J_{W},\deg (e_{I})\geq
i\}=\Lambda ^{*}W\cap \widehat{V_{i}}
\end{equation*}
where $J_{W}$ is the set of indices such that $W=span\{e_{i},i\in J_{W}\}.$
It follows that there is a unique notion of degree on $\Lambda ^{*}W$
associated to the original degree on $V$. 
\endproof%

In particular, this allows to define the degree of a subspace $W$ by setting 
$\deg _{V}(W)=\deg (f_{1}\wedge ...\wedge f_{k}),$ where $(f_{1},...,f_{k})$
is any basis of $W$ (observe that $\deg $ is really defined on the
projective space $\Bbb{P}(\Lambda ^{*}V)$). If $(e_{1},...,e_{k})$ is an
adapted basis for the filtration $(W\cap V_{i})_{i}$, then (\ref{degext})
yields 
\begin{equation}
\deg _{V}(W)=\deg (e_{1}\wedge ...\wedge e_{k})=\sum_{i=1}^{k}\deg
(e_{i})=\sum_{j\geq 1}\dim (W\cap V_{j}).  \label{deg(det)}
\end{equation}

Finally if $(\delta _{t})_{t>0}$ is the one-parameter group of endomorphisms
of $V$ defined by $\delta _{t}(e_{i})=t^{\deg (e_{i})}e_{i}$ for some
adapted basis $(e_{1},...,e_{n}).$ We easily check the following:

\begin{proposition}
\label{proplim}The $\delta _{t}$'s extend canonically to $\Lambda ^{*}V$ and
for any non-zero $\xi \in \Lambda ^{*}V$, $t^{-d}\delta _{t}(\xi )$ tends to
a non-zero limit in $\Lambda ^{*}V$ as $t\rightarrow 0$ if and only if $%
d=\deg (\xi ).$
\end{proposition}

\subsection{Homogeneous quasi-norms on filtered vector spaces}

Let $V$ be a real vector space with a filtration $V=V_{1}\supseteq
V_{2}\supseteq ...\supseteq V_{r}\supseteq V_{r+1}=\{0\}.$ We say that $%
(\delta _{t})_{t>0}$ is a \textit{one-parameter group of dilations}
associated to this filtration if the $\delta _{t}$'s are linear
automorphisms of $V$ such that, for every $i=1,...,r$, the eigenspace $%
m_{i}=\{x\in V,\delta _{t}(x)=t^{i}x\}$ is independent of $t$ ($t\neq 1$)
and is a (possibly trivial) supplementary vector subspace of $V_{i+1}$
inside $V_{i},$ i.e. 
\begin{equation}
V_{i}=m_{i}\oplus V_{i+1}.  \label{supplementary}
\end{equation}
The $r$ subspaces $m_{i}$'s are determined by $(\delta _{t})_{t>0}$ and
vice-versa any choice of $r$ subspaces $m_{i}$'s verifying (\ref
{supplementary}) determines a unique group of dilations. If $%
V=\bigoplus_{1\leq i\leq r}m_{i}$ and $V=\bigoplus_{1\leq i\leq
r}m_{i}^{\prime }$ are two choices of supplementary subspaces, then the
associated one-parameter groups of dilations satisfy the relation $\delta
_{t}^{\prime }=\phi \circ \delta _{t}\circ \phi ^{-1}$ where $\phi $ is the
coordinate change from the first direct sum to the second. Also we have
uniformly on bounded subsets of $V$: 
\begin{equation}
\lim_{t\rightarrow +\infty }\delta _{\frac{1}{t}}\circ \phi \circ \delta
_{t}=id.  \label{convphi}
\end{equation}

We now introduce the following definition.

\begin{definition}
\label{quasinorm}A continuous function $|\cdot |:V\rightarrow \Bbb{R}_{+}$
is called a homogeneous quasi-norm associated to the dilations $(\delta
_{t})_{t},$ or simply a \textbf{quasi-norm}, if it satisfies the following
properties:

$(i)$ $|x|=0\Leftrightarrow x=0.$

$(ii)$ $|\delta _{t}(x)|=t|x|$ for all $t>0.$
\end{definition}

Examples of quasi-norms are given by supremum quasi-norms of the type $%
|x|=\max_{p}\left\| \pi _{p}(x)\right\| _{p}^{1/p}$ where $\left\| \cdot
\right\| _{p}$ are ordinary norms on the vector space $m_{p}$ and $\pi _{p}$
is the projection on $m_{p}$ according to the decomposition $%
V=\bigoplus_{1\leq i\leq r}m_{i}$. More examples will be given below in the
context of nilpotent Lie groups.

Clearly, two quasi-norms associated to the same one-parameter group of
dilations are equivalent in the sense that $\frac{1}{c}\left| \cdot \right|
_{1}\leq \left| \cdot \right| _{2}\leq c\left| \cdot \right| _{1}$ for some
constant $c>0$. Furthermore, using (\ref{convphi}), we check the following:

\begin{proposition}
If $\left| \cdot \right| ^{\prime }$ is a quasi-norm that is homogeneous
with respect to the one-parameter group $(\delta _{t}^{\prime })_{t}$
associated to the decomposition $V=\bigoplus_{1\leq i\leq r}m_{i}^{\prime }$
then there is a unique quasi-norm $|\cdot |$ that is homogeneous with
respect to $(\delta _{t})_{t}$ such that 
\begin{equation*}
|x|-|x|^{\prime }=o(|x|)
\end{equation*}
as $|x|$ is large. In fact $|x|=|\phi (x)|^{\prime }.$
\end{proposition}

As will be observed below in Section \ref{Comp}, to every reasonable
left-invariant distance on a simply connected nilpotent Lie group $N$ is
associated a unique quasi-norm that is asymptotic to it. In particular,
every large ball for a left-invariant distance is well approximated by some
quasi-norm ball. This fact makes the volume computations needed in our main
theorem possible because, thanks to their scaling property, such
computations are easy in the case of quasi-norm balls.

\subsubsection{Invariance under restriction to a subspace or projection to a
quotient}

Let $|\cdot |$ be a homogeneous quasi-norm associated to some fixed
one-parameter group of dilations $(\delta _{t})_{t>0}$ and let $%
D_{t}:=\{x\in V,|x|\leq t\}$ be the corresponding quasi-norm ball. Let $W$
be a vector subspace of $V$ endowed with the induced filtration $(W\cap
V_{i})_{i}$ and let $(\delta _{t}^{0})_{t}$ be some one-parameter group of
dilations on $W$ with respect to that filtration. Although the restriction
of $|\cdot |$ to $W$ may not be a quasi-norm on $W,$ the following holds:

\begin{proposition}
There exists a unique homogeneous quasi-norm $|\cdot |_{0}$ on $W$ such that 
$|x|-|x|_{0}=o(|x|)$ if $x\in W$ and $|x|$ is large. In particular, if $%
D_{t}^{0}$ is the homogeneous quasi-norm ball for $|\cdot |_{0}$ on $W,$
then there exists $\varepsilon _{t}>0,$ $\varepsilon _{t}\rightarrow 0$ as $%
t\rightarrow +\infty ,$ such that 
\begin{equation*}
D_{t(1-\varepsilon _{t})}^{0}\subset D_{t}\cap W\subset D_{t(1+\varepsilon
_{t})}^{0}.
\end{equation*}
\end{proposition}

Similarly, we can consider the quotient vector space $V/W$ endowed with the
induced filtration $(V_{i}/V_{i}\cap W)_{i}$ and with some choice of a
one-parameter group of dilations $(\overline{\delta }_{t})_{t}$. Let $\pi
:V\rightarrow V/W$ be the canonical projection and let $|y|_{\pi }:=\inf
\{|x|,\pi (x)=y\}$. Although $|\cdot |_{\pi }$ may not be a quasi-norm on $%
V/W,$ the following holds:

\begin{proposition}
There exists a unique homogeneous quasi-norm $|\cdot |_{1}$ on $V/W$ such
that $|y|_{\pi }-|y|_{1}=o(|y|_{\pi })$ if $y\in V/W$ and $|y|_{\pi }$ is
large. In particular, if $D_{t}^{1}$ is the homogeneous quasi-norm ball for $%
|\cdot |_{1}$ on $V/W,$ then there exists $\varepsilon _{t}>0,$ $\varepsilon
_{t}\rightarrow 0$ as $t\rightarrow +\infty ,$ such that 
\begin{equation*}
D_{t(1-\varepsilon _{t})}^{1}\subset \pi (D_{t})\subset D_{t(1+\varepsilon
_{t})}^{1}.
\end{equation*}
\end{proposition}

We leave the proof of these propositions as an exercise.

\subsection{Volume growth\label{volgrowth}}

Let $N$ be a simply connected nilpotent Lie group and $(C^{p}(N))_{p=1,...,r}
$ its descending central series$.$ The integer $r$ is the nilpotency length,
that is the largest $r$ for which $C^{r}(N)$ is non trivial. We identify $N$
with its Lie algebra $\frak{n}$ via the exponential map, which is a
diffeomorphism. The Lie product is a polynomial function on $\frak{n}\times 
\frak{n}$ and this makes $N$ into a real algebraic group with a Zariski
topology. By Theorem 2.1 of \cite{Rag}, a closed subgroup of $N$ is
Zariski-dense if and only if it is co-compact. Furthermore, any 
closed subgroup $H$
of $N$ is contained in a unique minimal connected closed subgroup $%
\widetilde{H}$ of $N.$ The subgroup $\widetilde{H}$ is Zariski closed,
simply connected, and $\widetilde{H}/H$ is compact.

It was proved by Guivarc'h in \cite{Guiv} (and independently by Bass \cite
{Bass} in the special case of finitely generated nilpotent groups) that if $G
$ is a closed co-compact subgroup of $N$ and $U$ be a compact generating
neighborhood of the identity in $G$, then there are positive constants $C_{1}
$ and $C_{2}$ such that for any positive integer $n$%
\begin{equation}
C_{1}\cdot n^{d(N)}\leq vol_{G}(U^{n})\leq C_{2}\cdot n^{d(N)}  \label{ineq}
\end{equation}
where $d(N)$ is an integer called the \textit{homogeneous dimension} of $N$
and is given by the Bass-Guivarc'h formula: 
\begin{equation}
d(N)=\sum_{p\geq 1}\dim (C^{p}(N))=\sum_{p\geq 1}p\cdot \dim
(C^{p}(N)/C^{p+1}(N)).  \label{exp}
\end{equation}
For general $H$ as above we set $d(H)=d(\widetilde{H}).$ For instance, if $%
\Gamma $ is any nilpotent group generated by a finite symmetric set $S$
then, we can view it, modulo its finite torsion group, as a lattice in a
simply connected nilpotent Lie group $N$ according to a theorem of Malcev (%
\cite{Rag} chp. 2). Hence $d(\Gamma )=\sum_{p\geq 1}p\cdot rk(C^{p}(\Gamma
)/C^{p+1}(\Gamma ))$ and by (\ref{ineq}) the ball $S^{n}$ of radius $n$ in
the word metric defined $S$ has, up to multiplicative constants, $%
n^{d(\Gamma )}$ elements.

The estimate (\ref{ineq}) was later refined by Pansu who showed in \cite{Pan}
that $|S^{n}|/n^{d(\Gamma )}$ has a non-zero limit when $n\rightarrow
+\infty $. In fact, Pansu's argument can be adapted (see \cite{Bre2} for
details) to extend his result to all closed subgroups of $N,$ namely:

\begin{theorem}
\label{volumelimit}Let $G$ be a closed subgroup of $N$ and $U$ be a compact
generating neighborhood of the identity in $G$. Then there is a positive
constant $AsVol(U)>0$%
\begin{equation*}
\lim_{n\rightarrow +\infty }\frac{vol_{G}(U^{n})}{n^{d(G)}}=AsVol(U)\cdot 
\end{equation*}
\end{theorem}

Note that the homogeneous dimension $d(N)$ coincides the degree $\deg _{%
\frak{n}}(\frak{n})$ defined in (\ref{deg(det)}) where the filtration on $%
\frak{n}=Lie(N)$ is given by the central descending series. If $M$ is a
normal closed and connected subgroup of $N,$ then we check that $\deg _{%
\frak{n}}(\frak{m})=d(N)-d(N/M).$

\subsection{Polynomials, dilations and quasi-norms}

We recall here a few well-known facts about the analysis on nilpotent Lie
groups (see \cite{Goo}).

\subsubsection{Degree of an element and a polynomial\label{pol}}

Let $N$ be a simply connected nilpotent Lie group, which we identify to its
Lie algebra $\frak{n}$ via the exponential map $exp:\frak{n}\rightarrow N$.
We say that a map $P:N\rightarrow \Bbb{R}$ is polynomial if $x\mapsto P(\exp
(x))$ is a polynomial map on the real vector space $\frak{n}$. The central
descending series $(C^{k}(\frak{n}))_{k\geq 1}$ gives a canonical filtration
on $\frak{n}$ and hence induces a degree function on $\frak{n}$ (as in Prop. 
\ref{degprop}).

Let $(e_{i})_{i=1,...,n}$ be an \textit{adapted basis} of $\frak{n}$, namely
we assume that $C^{i}(\frak{n})=span\{e_{k}|\deg (e_{k})\geq i\}$ for all $%
i\geq 1.$ We define the \textit{degree} of a monomial $x^{\alpha
}=x_{1}^{\alpha _{1}}\cdot ...\cdot x_{n}^{\alpha _{n}}$ to be $d(\alpha
):=\alpha _{1}\deg (e_{1})+...+\alpha _{n}\deg (e_{n})$ for any multi-index $%
\alpha =(\alpha _{1},...,\alpha _{n})$, and the degree of an arbitrary
polynomial map to be the maximum degree of each of its monomials. This
definition is easily seen to be independent of the choice of the adapted
basis used to define it.

The coordinates of the product of two elements in the basis $%
(e_{i})_{i=1,..,n}$ are obtained from the Campbell-Hausdorff formula as
follows (see \cite{Goo} p. 14): 
\begin{equation}
(xy)_{i}=x_{i}+y_{i}+P_{i}(x,y)  \label{Campbell}
\end{equation}
where $P_{i}$ is a polynomial map on $\frak{n}\times \frak{n}$ of the
following special type: 
\begin{equation*}
P_{i}(x,y)=\sum C_{\alpha ,\beta }x^{\alpha }y^{\beta }
\end{equation*}
where $d(\alpha )+d(\beta )\leq \deg (e_{i})$ and $d(\alpha )\geq 1,d(\beta
)\geq 1$ and some constants $C_{\alpha ,\beta }.$

\subsubsection{Associated graded algebra, dilations\label{dilations}}

Let $(m_{p})_{p\geq 1}$ be a collection of supplementary subspaces on $\frak{%
n}$ as in (\ref{supplementary}) for $V_{i}=C^{i}(\frak{n).}$ If $x\in \frak{n%
},$ we write $x=\sum_{p\geq 1}\pi _{p}(x)$ where $\pi _{p}(x)$ is the linear
projection onto $m_{p}$. Let $(\delta _{t})_{t>0}$ be the one-parameter
group of dilations associated to the $m_{p}$'s. The dilations $\delta _{t}$
do not \textit{a priori} preserve the Lie bracket on $\frak{n}$. This is the
case if and only if 
\begin{equation}
\lbrack m_{p},m_{q}]\subseteq m_{p+q}  \label{gradation}
\end{equation}
for every $p$ and $q$. If (\ref{gradation}) holds, we say that the $%
(m_{p})_{p\geq 1}$ form a \textit{gradation} of the Lie algebra $\frak{n}$,
and that $\frak{n}$ is a \textit{graded} Lie algebra and $N$ is called a
Carnot group.

If (\ref{gradation}) does not hold, we can nevertheless consider a new Lie
algebra structure on the real vector space $\frak{n}$ by setting $%
[x,y]_{0}=\pi _{p+q}([x,y])$ if $x\in m_{p}$ and $y\in m_{q}$. This new
structure $\frak{n}_{0}$ is graded and the $(\delta _{t})_{t>0}$ are
automorphisms. We denote by $N_{0}$ the associated Lie group. In fact the
original Lie bracket $[x,y]$ on $\frak{n}$ can be deformed continuously to $%
[x,y]_{0}$ by setting $[x,y]_{t}=\delta _{t}([\delta _{\frac{1}{t}}x,\delta
_{\frac{1}{t}}y])$ and letting $t\rightarrow 0$.

On the other hand, the \textit{graded Lie algebra} associated with $\frak{n}$
is by definition $gr(\frak{n})=\bigoplus_{p\geq 1}C^{p}(\frak{n})/C^{p+1}(%
\frak{n}).$ endowed by the Lie bracket induced from the Lie bracket of $%
\frak{n}$. The quotient map $m_{p}\rightarrow C^{p}(\frak{n})/C^{p+1}(\frak{n%
})$ gives rise to a linear isomorphism between $\frak{n}$ and $gr(\frak{n})$%
, which is a Lie algebra isomorphism between $\frak{n}_{0}$ and $gr(\frak{n}%
).$ Hence graded Lie algebra structures induced by a choice of supplementary
subspaces $(m_{p})_{p\geq 1}$ as in (\ref{supplementary}) are all isomorphic
to $gr(\frak{n}).$

\subsubsection{Homogeneous quasi-norms on nilpotent Lie groups}

Let $|\cdot |$ be a homogeneous quasi-norm on $N$ with respect to some
direct sum decomposition given by $(m_{p})_{p}$'s$.$ Let $(e_{i})_{i}$ be an
adapted basis of $\frak{n}$.

\begin{proposition}
\label{propofqn}There is a constant $C>0$ such that

$(a)$ $|x_{i}|\leq C\cdot |x|^{\deg (e_{i})}$ if $x=x_{1}e_{1}+...+x_{n}e_{n}
$.

$(b)$ $|x^{-1}|\leq C\cdot |x|.$

$(c)$ $|xy|\leq C(|x|+|y|+1).$
\end{proposition}

It is straightforward to check $(a)$ and $(b)$. It can be a problem that the
constant in $(c)$ need not be $1$. In fact this is why we use the word
quasi-norm instead of just norm: we do not require the triangle inequality
axiom to hold. However the following lemma of Guivarc'h is often a good
enough remedy to this situation. Let $\left\| \cdot \right\| _{p}$ be an
arbitrary norm on the vector space $m_{p}$.

\begin{lemma}
\label{GuivLem}\label{jauge}(\cite{Guiv} lemme II.1) Up to rescaling each $%
\left\| \cdot \right\| _{p}$ into a proportional norm $\lambda _{p}\left\|
\cdot \right\| _{p}$ ($\lambda _{p}>0$) if necessary, the quasi-norm $%
|x|=\max_{p}\left\| \pi _{p}(x)\right\| _{p}^{1/p}$ satisfies $|xy|\leq
|x|+|y|+c$ for some constant $c>0$ and for all $x,y\in N$. Besides $N$ is
graded with respect to $(\delta _{t})_{t}$ if and only if $c=0$.
\end{lemma}

The proof is based on the Campbell-Hausdorff formula $(\ref{Campbell})$.
Lemma \ref{GuivLem} yields property $(c)$ above and also is the key step to
prove (\ref{ineq}).

\begin{example}
Note that one important class of quasi-norms consists of those of the form $%
|x|=d(e,x)$, where $d(x,y)$ is a Carnot-Carath\'{e}odory Finsler metric
induced on a graded nilpotent Lie group by some ordinary norm on the vector
subspace $m_{1}$.
\end{example}

\subsection{Nicely growing subsets\label{nicely}}

We define here \textit{nicely growing subsets}. These are essentially Folner
subsets with some extra properties that behave well under intersection with
a connected subgroup. Our main result, Theorem \ref{SecWeyl} below, will
hold for all such families of subsets.

Recall that $\deg _{N}$ is the degree function from Paragraph \ref{degreesub}
defined for all vector subspaces of $\frak{n}=Lie(N)$.

\begin{definition}
\label{nicegrowth}We say that a family of measurable subsets $(A_{t})_{t>0}$
of $N$ is \textbf{nicely growing} if it satisfies the following properties:

$(i)$ $(A_{t})_{t>0}$ increases and exhausts $N,$ i.e. $A_{t}\subseteq A_{s}$
if $t\leq s$ and $\bigcup_{t>0}A_{t}=N$.

$(ii)$ For every compact subset $K\subset N,$ there exists a positive
function $\varepsilon _{t}>0$ with $\varepsilon _{t}\rightarrow 0$ as $%
t\rightarrow +\infty $ such that $A_{t(1-\varepsilon _{t})}\subseteq
KA_{t}K\subseteq A_{t(1+\varepsilon _{t})}$ for all $t$ large enough.

$(iii)$ For any connected subgroup $M$ of $N,$ there exists a constant $%
C(M)>0$ such that 
\begin{equation*}
\lim_{t\rightarrow +\infty }\frac{vol_{M}(A_{t}\cap M)}{t^{\deg _{N}(M)}}%
=C(M).
\end{equation*}

$(iv)$ There is a constant $\alpha >1$ such that $A_{t}A_{t}^{-1}\subseteq
A_{\alpha t}$ for all $t>1$.

Additionally, if $\Gamma $ is a finitely generated torsion free nilpotent
group, then a family $(\Lambda _{t})_{t>0}$ of subsets of $\Gamma $ is said
to be \textbf{nicely growing} if there is a nicely growing family $%
(A_{t})_{t>0}$ of subsets of the Malcev closure of $\Gamma $ such that $%
\Gamma \cap A_{t(1-\varepsilon _{t})}\subseteq \Lambda _{t}\subseteq \Gamma
\cap A_{t(1+\varepsilon _{t})}$ for all $t>0$ and some positive function $%
\varepsilon _{t}>0$ with $\varepsilon _{t}\rightarrow 0$ as $t\rightarrow
+\infty .$
\end{definition}

It will be convenient to broaden this definition a little bit by allowing
different degree functions than $\deg _{N}$ in axiom $(iii).$ By a degree
function on $N$, we mean any map $\deg $ from the set of all connected
subgroups of $N$ to $\Bbb{N}$ which is non-decreasing in the sense that $%
M_{1}\subset M_{2}\Rightarrow \deg M_{1}\leq \deg M_{2}.$ Then we can speak
of a nicely growing family of subsets of $N$ \textit{relative to} the degree
function $\deg .$ By definition, such a family of subsets will satisfy all
four axioms except that axiom $(iii)$ will now hold with $\deg (M)$ in place
of $\deg _{N}(M).$ The notion of nicely growing sets is stable under
intersection with a connected subgroup say $M$, but then the degree function
remains $\deg _{N}.$

Let us recall the statement of Proposition \ref{propdeb}, which is our main
goal here because it provides us with the many examples of nicely growing
subsets.

\begin{proposition}
\label{ThmDiscs}\label{evenseq}For any quasi-norm on $N,$ the balls $%
(D_{t})_{t>0}$, $D_{t}=\{|x|\leq t\},$ form a nicely growing family of
subsets of $N$.
\end{proposition}

\proof%
Property $(i)$ is obvious and property $(iv)$ follows from properties $(b)$
and $(c)$ of Prop. \ref{propofqn}. As for property $(ii),$ it is a
consequence of the following more general fact.

\begin{lemma}
\label{NormCont}Let $|\cdot |$ be a quasi-norm on $N$. Then for every $%
\varepsilon >0$ there exists $\delta >0$ and $C>0$ such that if $|y|>C$ and $%
|x|\leq \delta |y|$ we have 
\begin{equation*}
\left| |xy|-|y|\right| \leq \varepsilon |y|.
\end{equation*}
\end{lemma}

\proof%
Let $|\cdot |_{0}$ be the quasi-norm obtained from Guivarch's lemma (Lemma 
\ref{GuivLem}). Since any two quasi-norms are equivalent, $|x|_{0}=o(|y|_{0})
$ and $|x|=o(|y|)$ are equivalent conditions. Using Prop. \ref{propofqn} $(b)
$ and Lemma \ref{GuivLem}, we see that if $|x|_{0}=o(|y|_{0})$ then 
\begin{equation}
\left| |xy|_{0}-|y|_{0}\right| =o(|y|_{0}).  \label{fornought}
\end{equation}
Now write $\left| |xy|-|y|\right| =\left| \left| \delta _{\frac{1}{|xy|_{0}}%
}(xy)\right| |xy|_{0}-\left| \delta _{\frac{1}{|y|_{0}}}(y)\right|
|y|_{0}\right| .$ By (\ref{fornought}), we get $\left| |xy|-|y|\right|
=\left| \left| \delta _{\frac{1}{|y|_{0}}}(xy)\right| |xy|_{0}-\left| \delta
_{\frac{1}{|y|_{0}}}(y)\right| |y|_{0}\right| +o(|y|).$ So it remains to
show that $\delta _{\frac{1}{|y|_{0}}}(xy)\delta _{\frac{1}{|y|_{0}}}(y^{-1})
$ tends to $1$ as $|y|$ tends to infinity and $|x|=o(|y|).$ By compactness,
we may even assume that $\delta _{\frac{1}{|y|_{0}}}(y)=z$ is fixed and,
posing $t=|y|_{0},$ we have reduced to showing that $\delta _{\frac{1}{t}%
}(x\delta _{t}(z))\rightarrow z$ as $t\rightarrow +\infty $ and $|x|=o(t).$
To do this we use formula (\ref{Campbell}) for the product in coordinates 
\begin{equation*}
\left[ \delta _{\frac{1}{t}}(x\delta _{t}(z))\right] _{i}=\frac{x_{i}}{%
t^{\deg (e_{i})}}+z_{i}+\sum C_{\alpha ,\beta }x^{\alpha }\frac{z^{\beta }}{%
t^{\deg (e_{i})-d(\beta )}}
\end{equation*}
with the additional constraints $d(\alpha )+d(\beta )\leq \deg (e_{i})$ and $%
d(\alpha )\geq 1$, $d(\beta )\geq 1.$ The condition $|x|=o(t)$ means that $%
x_{i}=o(t^{\deg (e_{i})})$ for any index $i$. As $t$ tends to infinity, we
indeed obtain the convergence of the above expression towards $z_{i}.$ 
\endproof%

Observe that the analogous result holds when $xy$ is changed into $yx$ ($%
x\mapsto |x^{-1}|$ is another equivalent quasi-norm).

\medskip We now turn to the proof of property $(iii)$. This is where we will
need the discussion on filtrations from the previous sections. Let us denote
by $\frak{m}=Lie(M)$ the Lie subalgebra of $\frak{n}=Lie(N)$ corresponding
to $M$. Let $d=\dim (M)$ and $(f_{1},...,f_{d})$ a basis for the vector
space $\frak{m}$.

The orthogonal linear transformations (for some Euclidean norm $\left\|
\cdot \right\| $ on $\frak{n}$) act transitively on $Gr_{d}(\frak{n}),$ the
Grassmannian variety of $d$-dimensional linear subspaces of $\frak{n}$.
Hence there is some orthogonal map $o_{t}$ such that $o_{t}^{-1}\delta _{%
\frac{1}{t}}$ fixes $\frak{m}$. Let $vol_{M}$ be a Haar measure on $M,$
which we identify with Lebesgue measure on $\frak{m}.$ Let $\alpha _{t}^{-1}$
be the abolute value of the determinant of the endomorphism induced on $%
\frak{m}$ by $o_{t}^{-1}\delta _{\frac{1}{t}}.$ We have
\begin{equation*}
\alpha _{t}^{-1}=\frac{\left\| \delta _{\frac{1}{t}}f_{1}\wedge ...\wedge
\delta _{\frac{1}{t}}f_{d}\right\| }{\left\| f_{1}\wedge ...\wedge
f_{d}\right\| }
\end{equation*}
We can estimate the behavior of $\alpha _{t}$ when $t\rightarrow +\infty $.
As follows from (\ref{deg(det)}) above, 
\begin{equation*}
\deg (f_{1}\wedge ...\wedge f_{d})=\sum_{i\geq 1}\dim (\frak{m}\cap C^{i}(%
\frak{n}))=\deg _{N}(M)
\end{equation*}
Hence by Prop. \ref{proplim} there is a non-zero $\xi \in \Lambda ^{d}\frak{n%
}$ such that 
\begin{equation}
\lim_{t\rightarrow +\infty }t^{\deg _{N}(M)}\delta _{\frac{1}{t}}f_{1}\wedge
...\wedge \delta _{\frac{1}{t}}f_{d}=\xi   \label{convgrass}
\end{equation}
Then we can define 
\begin{equation}
c_{M}:=\lim_{t\rightarrow +\infty }\frac{\alpha _{t}}{t^{\deg _{N}(M)}}=%
\frac{\left\| f_{1}\wedge ...\wedge f_{d}\right\| }{\left\| \xi \right\| }>0
\label{cm}
\end{equation}
By (\ref{convgrass}) the subspaces $\delta _{\frac{1}{t}}\frak{m}$ converge
to a limit subspace $\frak{m}_{\infty }$ in the Grassmannian variety. So we
could choose $o_{t}$ so that $o_{t}$ converges to some $o$ as $t\rightarrow
+\infty .$ Then $\frak{m}_{\infty }=o\frak{m}$. Also observe that $\frak{m}%
_{\infty }$ is invariant under the full one-parameter group of dilations $%
(\delta _{t})_{t}$.

\begin{lemma}
\label{Lemm1}The subspace $\frak{m}_{\infty }$ is a Lie subalgebra of $\frak{
n}_{0}.$
\end{lemma}

\proof
We need to show that if $x,y\in \frak{m}_{\infty }$, then $x\cdot y\in \frak{%
m}_{\infty }$, where $x \cdot y$ is the product in $\frak{n}_{0}$. Let $x_{t},y_{t}$ in $\frak{m}$ be such that $%
x=\lim_{t\rightarrow +\infty }\delta _{\frac{1}{t}}(x_{t})$ and $%
y=\lim_{t\rightarrow +\infty }\delta _{\frac{1}{t}}(y_{t}).$ Then by (\ref
{Campbell}) we have 
\begin{equation*}
\left[ \delta _{\frac{1}{t}}(x_{t}y_{t})-\delta _{\frac{1}{t}%
}(x_{t})\cdot \delta _{\frac{1}{t}}(y_{t})\right] _{i}=\frac{1}{t^{\deg
(e_{i})}}\sum_{d_{\alpha}+d_{\beta} < \deg(e_{i})} C_{\alpha ,\beta }x_{t}^{\alpha }y_{t}^{\beta }.
\end{equation*}
By Prop. \ref{propofqn} (a) this expression is a $O(\frac{1}{t}).$ Hence $x\cdot y=\lim_{t\rightarrow +\infty }\delta _{\frac{1}{t}}(x_{t}y_{t})$
, i.e. $x\cdot y\in \frak{m}_{\infty }.$ 
\endproof

Let $vol_{M_{\infty }}$ be a Haar measure on $M_{\infty }$, again identified
with Lebesgue measure on $\frak{m}_{\infty }.$ Fixing the Haar measure on $N,
$ the choice of a Euclidean norm on $\frak{n}$ specifies a normalization for 
$vol_{M}$ and $vol_{M_{\infty }}.$ Then we have:

\begin{lemma}
\label{Lemm2}We have the following weak convergence of measures:
\begin{equation*}
\lim_{t\rightarrow +\infty }\frac{1}{t^{\deg _{N}(\frak{m})}}\left( \delta _{%
\frac{1}{t}}\right) _{*}vol_{M}=c_{M}\cdot vol_{M_{\infty }}
\end{equation*}
\end{lemma}

\proof%
By (\ref{cm}) we have  
\begin{equation*}
\lim_{t\rightarrow +\infty }\frac{1}{t^{\deg _{N}(\frak{m})}}\left( \delta _{%
\frac{1}{t}}\right) _{*}vol_{M}=\lim_{t\rightarrow +\infty }\frac{\alpha _{t}%
}{t^{\deg _{N}(\frak{m})}}\left( o_{t}\right) _{*}vol_{M}=c_{M}\cdot
vol_{M_{\infty }}
\end{equation*}
\endproof%

By Lemma \ref{Lemm2} applied to $D_{1},$ we get the desired convergence,
namely 
\begin{equation*}
\lim_{t\rightarrow +\infty }\frac{vol_{M}(M\cap D_{t})}{t^{\deg _{N}(\frak{m}%
)}}=c_{M}\cdot vol_{M_{\infty }}(D_{1})
\end{equation*}
after we check that $vol_{M_{\infty }}(\partial D_{1})=0.$ This is clear
however since, using the invariance of $M_{\infty }$ under $(\delta
_{t})_{t>0},$ the function $vol_{M_{\infty }}(D_{t})=|Jac((\delta
_{t})_{|M_{\infty }})|\cdot vol_{M_{\infty }}(D_{1})$ is continuous as a
function of $t.$ This ends the proof of Proposition \ref{ThmDiscs}. 
\endproof%

\subsection{Other examples of nicely growing subsets}

Other types of balls give rise to nicely growing subsets besides quasi-norm
balls. We show here two more examples: balls obtained by considering
exponential coordinates of the second kind, and more importantly $\rho $%
-balls for any ``reasonable'' left-invariant distance $\rho $ on $N.$

\subsubsection{Balls in privileged coordinates}

Let $(e_{i})_{1\leq i\leq n}$ be an adapted basis as in Paragraph \ref{pol},
and $(x_{i})_{i}$ be the associated exponential coordinates of the first
kind, i.e. $x=x_{1}e_{1}+...+x_{n}e_{n}$. For $i=1,...,n$, let $P_{i}$ be a
polynomial map on $\frak{n}$ of total (homogeneous) degree $\leq d_{i}=\deg
(e_{i})$, with $P_{i}(0)=0.$ The map $P_{i}$ can be split in two parts, $%
P_{i}(x)=L_{i}(x)+M_{i}(x)$ where $L_{i}(x)$ is a linear form on $\frak{n}$
depending only on those coordinates $x_{j}$ such that $d_{j}=d_{i}$ and
where $M_{i}(x)$ is a polynomial depending only on those $x_{j}$'s with $%
d_{j}\leq d_{i}-1$. Assume further that the $L_{i}$'s, $i=1,...,n$, are
linearly independent. Then, following \cite{Bel}, we call \textit{privileged
coordinates }any choice of coordinates on $\frak{n}$ that is obtained from
the $x_{i}$'s by a coordinate change $\phi :x\mapsto x^{\prime }$ of the
form 
\begin{equation}
x_{i}^{\prime }=P_{i}(x)  \label{privilege}
\end{equation}
For instance, writing $x$ in the associated exponential coordinates of the
second kind, i.e. 
\begin{equation*}
x=\exp (x_{1}^{\prime }e_{1})\cdot ...\cdot \exp (x_{n}^{\prime }e_{n})
\end{equation*}
it is easy to check that (\ref{privilege}) holds and we even have $%
L_{i}(x)=x_{i}.$ We have:

\begin{proposition}
Let $\left| \cdot \right| $ be a homogeneous quasi-norm on $N$ and $\phi
:x\mapsto x^{\prime }$ a privileged coordinate change. Then there exists a
unique homogeneous quasi-norm $\left| \cdot \right| ^{\prime }$ such that $%
|\phi (x)|=|x|^{\prime }+o(|x|^{\prime }).$ In particular the balls $\left\{
|\phi (x)|\leq t\right\} $ form a nicely growing family if $t$ tends to $%
+\infty .$
\end{proposition}

\proof%
Let us write $\phi (x)=\psi (x)+\eta (x)$ where $(\psi (x))_{i}$ is the
component of $P_{i}(x)$ that is homogeneous of homogeneous degree $d_{i}$.
Then clearly, the map $x\mapsto |\psi (x)|$ is a homogeneous quasi-norm
(property $(i)$ of Def. \ref{quasinorm} follows from the linear independence
of the $L_{i}$'s). Therefore, after composing by $\psi ^{-1}$, we may assume
that $\psi (x)=x,$ i.e. $\phi (x)=x+\eta (x)$ where each $(\eta (x))_{i}$ is
a polynomial map of degree $\leq d_{i}-1$. We wish to show that $|\phi
(x)|=|x|+o(|x|).$ Since $|\eta (x)|=o(|x|)$ the proof is a straightforward
copy of that of Lemma \ref{NormCont}.

Given $\varepsilon >0$, we thus get $D_{t(1-\varepsilon )}\subset \left\{
|\phi (x)|\leq t\right\} \subset D_{t(1+\varepsilon )}$, for all large $t,$
where $D_{t}=\{|x|\leq t\}.$ By Prop. \ref{ThmDiscs}, this clearly implies
all of the four defining properties for nicely growing subsets.%
\endproof%

\subsubsection{Balls for left-invariant metrics\label{Comp}}

Following \cite{Bre2}, we say that a distance function $\rho $ on $N$ is a 
\textit{periodic metric} if it is left-invariant under some co-compact
subgroup of $N$ and is asymptotically geodesic, namely for every $%
\varepsilon >0$ there is $s>0$ such that $\forall x,y\in N,$ one can find
points $x_{1}=x,$ $x_{2},...,x_{n}=y$ in $N$ such that $\rho (x,y)\geq
(1-\varepsilon )\sum_{i=1}^{n-1}\rho (x_{i},x_{i+1})$ and $\rho
(x_{i},x_{i+1})\leq s$ for each $i.$ This means that we require a kind of
weak existence of geodesic axiom. Examples of such metrics are given by
left-invariant Riemmanian or sub-Riemmanian metrics on $N$, and also by word
metrics induced by a compact generating set of $N$ (see \cite{Bre2}).

In \cite{Pan}, Pansu associates to any such $\rho $ a Carnot-Carath\'{e}%
odory metric $d_{\rho }$ on $N$ in the following way. Let $(\delta _{t})_{t}$
be a one-parameter group of (linear) dilations on $\frak{n}$ and $\frak{n}%
=\oplus _{p}m_{p}$ the eigenspace decomposition, with $\pi _{1}$ the
projection to $m_{1}$. This yields a graded structure $\frak{n}_{0}$ on $%
\frak{n}$ as defined in \ref{dilations}. Let $E_{s}$ is the closed convex
hull of all $\pi _{1}(x)/\rho (e,x)$ with $x\in N$ and $\rho (e,x)>s$ and $%
E=\bigcap_{s>0}E_{s}.$ Pansu shows that $E$ is a compact symmetric subset of 
$m_{1}$ with non-empty interior. Hence it is the unit ball of some norm $%
\left\| \cdot \right\| _{\rho }$ on $m_{1}.$ In turn, this norm defines a
Carnot-Carath\'{e}odory Finsler distance $d_{\rho }$ on $N_{0}$ by setting $%
d_{\rho }(x,y)=\inf \{L(\gamma ),\gamma $ horizontal path from $x$ to $y\}$,
where horizontal means almost everywhere tangent to a $\frak{n}_{0}$-left
translate of $m_{1}$ and $L(\gamma )$ is the length of $\gamma $ measured
according to $\left\| \cdot \right\| _{\rho }.$ Pansu's main result reads
(see \cite{Pan}, or \cite{Bre2} for a proof):.

\begin{theorem}
(\cite{Pan}) \label{Pansu1}For any periodic metric $\rho $ on $N,$%
\begin{equation}
\lim_{x\rightarrow \infty }\frac{\rho (e,x)}{d_{\rho }(e,x)}=1  \label{asym}
\end{equation}
\end{theorem}

\begin{corollary}
\label{Panpan}The balls $B_{\rho }(t)=\{x\in N,\rho (e,x)\leq t\}$ form a
family of nicely growing subsets of $N$.
\end{corollary}

\proof%
Clearly $x\mapsto d_{\rho }(e,x)$ is a homogeneous quasi-norm, so the
Corollary follows from Prop. \ref{ThmDiscs}. 
\endproof%

When proving Theorem \ref{main} we will apply this result to the case when $%
\rho (x,y)=d(\gamma _{x},\gamma _{y})$, where $d$ is any word metric on a
lattice $\Gamma $ in $N$, and $x\in \gamma _{x}F$ for some fixed compact
fundamental domain $F$ of $\Gamma $ in $N$.

\section{Equidistribution\label{SecWeyl}}

In this section, we prove Theorem \ref{main} via the following version of it:

\begin{theorem}
\label{eq}(Dense subgroups are equidistributed) Let $G$ be a closed subgroup
of a simply connected nilpotent Lie group and $\Gamma $ a finitely generated
torsion-free nilpotent group. Let $\phi :\Gamma \rightarrow G$ be a
homomorphism with dense image$.$ Suppose that $(\Lambda _{t})_{t>0}$ is a
nicely growing family of subsets of $\Gamma $ (see Definition \ref
{nicegrowth})$.$ Then there is a positive constant $C_{1}>0$ depending on $%
(\Lambda _{t})_{t>0}$ and on the choice of a Haar measure $vol_{G}$ on $G,$
such that for any bounded Borel subset $B\subset G$ with negligible boundary
we have 
\begin{equation}
\frac{\#\left\{ \gamma \in \Gamma ,\gamma \in \Lambda _{T},\phi (\gamma )\in
B\right\} }{T^{d(\Gamma )-d(G)}}\underset{T\rightarrow +\infty }{\rightarrow 
}C_{1}\cdot vol_{G}(B)  \label{est}
\end{equation}
where $d(\Gamma )$ and $d(G)$ are the integers defined in \ref{volgrowth}.
\end{theorem}

Recall that by Malcev's theory (see \cite{Rag} Theorem 2.18), every finitely
generated torsion-free nilpotent group embeds in a simply connected
nilpotent Lie group, its Malcev closure$.$ Let $N$ be the Malcev closure of $%
\Gamma $. Let $(D_{t})_{t>0}$ be the family of nicely growing subsets of $N$
such that $\Gamma \cap D_{t(1-\varepsilon _{t})}\subseteq \Lambda
_{t}\subseteq \Gamma \cap D_{t(1+\varepsilon _{t})}$. First observe that it
is enough to prove the theorem for sets $\Lambda _{t}$ of the form $\Lambda
_{t}=\Gamma \cap D_{t}$.

Before starting the proof of Theorem \ref{eq}, let us briefly explain how
one can also reduce to the case when $G$ is connected. Recall that since $G$
is closed in a simply connected nilpotent Lie group, it is co-compact in the
simply connected subgroup $\widetilde{G}$ (its Zariski-closure) defined in 
\ref{volgrowth}. By Malcev's rigidity (see \cite{Rag} Theorem 2.11) $\phi $
extends to an epimorphism $\phi :N\rightarrow \widetilde{G}$, which gives
rise to an isomorphism $N/N_{0}\overset{\simeq }{\rightarrow }\widetilde{G}%
/G^{\circ },$ where $N_{0}=\widetilde{\Gamma ^{\circ }}$ is normal in $N$
and contains $M=\ker \phi ,$ and $\Gamma ^{\circ }=\Gamma \cap \phi
^{-1}(G^{\circ })$. For every bounded Borel subset of $G$, the set $\left\{
\gamma \in \Gamma ,\gamma \in D_{t},\phi (\gamma )\in B\right\} $ can be
split into a finite number of translates of $\left\{ \gamma \in \Gamma
^{\circ },\gamma \in \gamma _{i}D_{t},\phi (\gamma )\in B_{i}\right\} $
where $\phi (\gamma _{i}^{-1})B_{i}\subset B$ and $\gamma _{i}\in \Gamma $
and $B_{i}\subset G^{\circ }.$ Since $G^{\circ }$ is simply connected and $%
\Gamma ^{\circ }$ dense in it, we may just as well work with these groups.
However $D_{t}\cap N_{0}$ is a nicely growing family of subsets of $N_{0}$
only relative to the degree function $\deg _{N}$ and not relative to $\deg
_{N_{0}}$ (see the remarks below Definition \ref{nicegrowth}). Nevertheless,
we show below that if $D_{t}$ is a nicely growing family of subsets of $N$
with respect to an arbitrary degree function, then under the hypothesis and
notation of Theorem \ref{eq} 
\begin{equation}
\frac{\#\left\{ \gamma \in \Gamma ,\gamma \in D_{t},\phi (\gamma )\in
B\right\} }{vol_{M}(M\cap D_{t})}\underset{t\rightarrow +\infty }{%
\rightarrow }C_{1}\cdot vol_{G}(B)  \label{convvol}
\end{equation}
We will thus assume below that $G$ is connected and simply connected.

\subsection{Unique ergodicity and counting}

By Malcev's rigidity (\cite{Rag} Theorem 2.11), $\phi $ extends to an
epimorphism $\phi :N\rightarrow G$. Let $M=\ker \phi $. We want to find the
asymptotics of the number of points of $\Gamma \cap D_{t}$ which lie in $%
\phi ^{-1}(B),\;$the inverse image of the bounded Borel subset $B$ by the
map $\phi $. Hence we are dealing with a counting problem, which we will
treat via ergodic theory. The use of ergodic theory to solve counting
problems is now standard (see for instance \cite{EMS}, \cite{EMM} and \cite
{Bab} for a survey of these techniques) and what we are going to present
here is yet another illustration of these ideas.

It is a fairly general principle in ergodic theory that the ergodic
properties of the action of a closed subgroup $H_{1}$ of a group $H$ on the
homogeneous space $H/H_{2}$ can be deduced from the ergodic properties of
the action of the closed subgroup $H_{2}$ on $H/H_{1}$ and vice-versa. Here
we will deduce the equidistribution of $\Gamma $ in $G\simeq N/M$ from the
equidistribution of an $M$-orbit on the nilmanifold $N/\Gamma $. In order to
do this, we first recall the following well-known theorem (see \cite{Sta}
Theorems 3.6 and 3.8, and also \cite{AuB} Lemma 5.1):

\begin{theorem}
\label{Green}(Unique ergodicity criterion for nilflows) Let $\Gamma $ be a
co-compact lattice in a simply connected nilpotent Lie group $N$ and let $M$
be a closed subgroup of $N$. The following are equivalent:

$(i)$ The subset $M\Gamma $ is dense in $N.$

$(ii)$ $M$ acts ergodically on $N/\Gamma $.

$(iii)$ The $M$-action on $N/\Gamma $ is uniquely ergodic.
\end{theorem}

Since $\Gamma $ is dense in $G\simeq N/M,$ Theorem \ref{Green} implies that
the action of $M$ on $N/\Gamma $ is uniquely ergodic, i.e. that the
normalized Haar measure $\nu $ on $N/\Gamma $ is the only $M$-invariant
probability measure on $N/\Gamma $. In order to translate the counting
problem into an equidistribution question, we introduce the following
counting function for $x\in N,$ 
\begin{equation*}
F_{t}^{B}(x)=\#\left\{ \gamma \in \Gamma ,x\gamma \in D_{t},\phi (x\gamma
)\in B\right\}
\end{equation*}
Note that $F_{t}^{B}$ is $\Gamma $-invariant on the right hand side, hence
it really defines a measurable function on the nilmanifold $N/\Gamma $. Note
further that the quantity we are interested in is precisely $%
F_{t}^{B}(e)=\#\left\{ \gamma \in \Gamma ,\gamma \in D_{t},\phi (\gamma )\in
B\right\} ,$ and our goal (i.e. (\ref{convvol})) is to prove the convergence
of $F_{t}^{B}(e)/V_{t}$, where $V_{t}=vol_{M}(M\cap D_{t})$. The main step
is to prove weak convergence ((\ref{conv12}) below) of the functions $%
F_{t}^{B}(x)/V_{t}$.

\subsection{Proof of Theorems \ref{eq} and \ref{main}.}

Since the Haar measure $\nu $ on $N/\Gamma $ is already normalized, the
choice of a Haar measure on $G,$ denoted by $vol_{G},$ determines uniquely a
Haar measure on the kernel $M,$ which we denote by $vol_{M}$. For every $%
g\in N$ we denote by $\nu _{t}^{g}$ the image under $\pi :N\rightarrow
N/\Gamma $ of the uniform probability measure supported on $g^{-1}D_{t}\cap
M $, i.e. 
\begin{equation}
\nu _{t}^{g}=\pi _{*}\left( \frac{1_{D_{t}\cap gM}(gy)}{vol_{M}(g^{-1}D_{t}%
\cap M)}vol_{M}(dy)\right)  \label{defnutg}
\end{equation}
We let $\psi $ be a continuous function on $N/\Gamma $ and we consider the
scalar product 
\begin{equation*}
\left\langle F_{t}^{B},\psi \right\rangle =\int_{N/\Gamma }\sum_{\gamma \in
\Gamma }1_{x\gamma \in D_{t}}1_{\phi (x\gamma )\in B}\psi (\overline{x})\nu
(d\overline{x})=\int_{N}1_{D_{t}\cap \phi ^{-1}(B)}(x)\psi (\overline{x}%
)vol_{N}(dx)
\end{equation*}
Decomposing the Haar measure on $N$ along the fibers of the projection $\phi
:N\rightarrow N/M,$ we obtain, 
\begin{eqnarray}
\left\langle F_{t}^{B},\psi \right\rangle
&=&\int_{G}1_{B}(g)\int_{M}1_{D_{t}}(gy)\psi (\overline{gy}%
)vol_{M}(dy)vol_{G}(dg)  \label{Ft} \\
&=&\int_{B}vol_{M}(g^{-1}D_{t}\cap M)\left( \int_{N/\Gamma }\psi (gz)\nu
_{t}^{g}(dz)\right) vol_{G}(dg)  \notag
\end{eqnarray}

In order to go further, we need the following proposition, which is the
consequence of the unique ergodicity of the $M$-action on $N/\Gamma $.

\begin{proposition}
\label{uniquergod}The following weak convergence of probability measures on $%
N/\Gamma $ holds uniformly when $g$ varies in compact subsets of $N.$ 
\begin{equation*}
\lim_{t\rightarrow +\infty }\nu _{t}^{g}=\nu 
\end{equation*}
\end{proposition}

\proof%
Let $\nu _{\infty }$ be a weak limit of $\nu _{t}^{g}$ as $t\rightarrow
+\infty $ and $g$ converges to some element in $N$. Since the $M$-action is
uniquely ergodic on $N/\Gamma $ by Theorem \ref{Green}, it is enough to show
that $\nu _{\infty }$ is invariant under $M.$ Let $\mu _{t}^{g}$ be the
probability measure on $N$ such that $\nu _{t}^{g}=\pi _{*}(\mu _{t}^{g})$
as defined in (\ref{defnutg}). The map $\pi _{*}:\mathcal{P}(N)\rightarrow 
\mathcal{P}(N/\Gamma )$ between spaces of probability measures is an $N$%
-equivariant contraction for the total variation norm, hence to show that $%
\nu _{\infty }$ is invariant under $M,$ it is enough to prove the following
lemma.

\begin{lemma}
\label{totalvariation}For any $h\in M,$ the following convergence holds
uniformly in $g$ as $g$ varies in compact subsets of $N$%
\begin{equation*}
\lim_{t\rightarrow +\infty }\left\| \delta _{h}*\mu _{t}^{g}-\mu
_{t}^{g}\right\| =0
\end{equation*}
\end{lemma}

\proof%
Since $h\in M$, the measures $\mu _{t}^{g}$ and $\delta _{h}*\mu _{t}^{g}$
are supported on $M$ and are absolutely continuous with respect to the Haar
measure $vol_{M}$. Hence the total variation norm is simply the $\Bbb{L}^{1}$
norm. So 
\begin{equation}
\left\| \delta _{h}*\mu _{t}^{g}-\mu _{t}^{g}\right\| =\frac{vol_{M}\left(
M\cap (h^{-1}g^{-1}D_{t}\Delta g^{-1}D_{t})\right) }{vol_{M}(M\cap
g^{-1}D_{t})}  \label{norm1}
\end{equation}
where $\Delta $ is the symmetric difference operator. Note that combining
both properties $(ii)$ and $(iii)$ of Definition \ref{nicegrowth}, the
following convergence holds uniformly when $g$ varies in compact subsets of $%
N$. 
\begin{equation}
\lim_{t\rightarrow +\infty }\frac{vol_{M}(g^{-1}D_{t}\cap M)}{t^{\deg (M)}}%
=C(M)  \label{norm2}
\end{equation}
Similarly, by property $(ii),$ if $g$ lies in a compact set, there will be
some positive function $\varepsilon _{t}>0$ with $\varepsilon
_{t}\rightarrow 0$ such that we can write $h^{-1}g^{-1}D_{t}\Delta
g^{-1}D_{t}\subseteq D_{t+t\varepsilon _{t}}\backslash D_{t-t\varepsilon
_{t}},$ But by property $(iii)$ we can conclude that 
\begin{equation}
\lim_{t\rightarrow +\infty }\frac{vol_{M}\left( (D_{t+t\varepsilon
_{t}}\backslash D_{t-t\varepsilon _{t}})\cap M\right) }{vol_{M}(M\cap D_{t})}%
=0  \label{norm4}
\end{equation}
Combining (\ref{norm1}) with (\ref{norm2}) and (\ref{norm4}) we are done. 
\endproof%

\endproof%

Let us resume the proof of Theorem \ref{eq}. Since $B$ is bounded and $g\in B
$ in the integral (\ref{Ft}), when $t$ tends to $+\infty $ we obtain the
weak convergence: 
\begin{equation}
\lim_{t\rightarrow +\infty }\frac{\left\langle F_{t}^{B},\psi \right\rangle 
}{vol_{M}(D_{t}\cap M)}=vol_{G}(B)\cdot \int_{N/\Gamma }\psi (z)d\nu (z)
\label{conv12}
\end{equation}
In order to get the desired asymptotics for $F_{t}^{B}(e)$, we need to
compare it to $\left\langle F_{t}^{B},\psi \right\rangle $ where $\psi $ is
chosen to better and better approximate the Dirac distribution $\delta _{e}$%
. For every sufficiently small neighborhood of the identity $U$ in $N,$
which is homeomorphic to $U$ via the covering $\pi $, we may consider a bump
function $\psi $ supported on $\pi (U)$ (i.e. a non-negative continuous
function with total sum equal to $1$). Then $F_{t}^{B}(e)\leq
F_{t+o(t)}^{\phi (U)B}(g)$ for any $g\in U$ by property $(ii)$ of nicely
growing subsets. Note also that $vol_{G}(\phi (U)B)$ gets closer and closer
to $vol_{G}(B)$ as $U$ tends to the identity. It follows that $%
F_{t}^{B}(e)\leq \left\langle F_{t+o(t)}^{\phi (U)B},\psi \right\rangle $
and from property $(iii)$ of nicely growing subsets, applying (\ref{conv12}%
), we obtain that 
\begin{equation}
\underset{t\rightarrow +\infty }{\overline{\lim }}\frac{F_{t}^{B}(e)}{%
vol_{M}(D_{t}\cap M)}\leq vol_{G}(B)  \label{123up}
\end{equation}

The lower bound is only slightly more delicate. We have $F_{t}^{B}(e)\geq
F_{t}^{\overset{\circ }{B}}(e)$ where $\overset{\circ }{B}$ is the interior
of $B.$ Since $B$ has negligible boundary, for any sufficiently small
neighborhood $V$ of the identity in $G,$ there exists an open subset $B_{V}$
of $B$ such that $VB_{V}\subset B$ and $vol_{G}(B_{V})\rightarrow vol_{G}(B)$
as $V$ narrows to the identity.

Now let $U$ be a neighborhood of the identity in $N$ so small that $\phi
(U^{-1})\subset V$. Then for any $g\in U,$ $F_{t}^{B}(e)\geq
F_{t-o(t)}^{B_{V}}(g)$ for any $g\in U$. It follows that $F_{t}^{B}(e)\geq
\left\langle F_{t-o(t)}^{B_{V}},\psi \right\rangle $ for any bump function
supported on $\pi (U)$. Again from property $(iii)$ of nicely growing
subsets, applying (\ref{conv12}), we obtain 
\begin{equation*}
\underset{t\rightarrow +\infty }{\underline{\lim }}\frac{F_{t}^{B}(e)}{%
vol_{M}(D_{t}\cap M)}\geq vol_{G}(B).
\end{equation*}
We have proved 
\begin{equation*}
\underset{t\rightarrow +\infty }{\lim }\frac{F_{t}^{B}(e)}{vol_{M}(D_{t}\cap
M)}=vol_{G}(B)
\end{equation*}
which, together with the volume estimate given by property $(iii)$ of
Definition \ref{nicegrowth}, ends the proof of Theorem \ref{eq}.%
\endproof%

\proof[Proof of Theorem \ref{main}]%
If $T$ is the torsion subgroup of $\Gamma $, the map $\phi $ factors through 
$T$ to a quotient map $\overline{\phi }:\Gamma /T\rightarrow G.$ Since $T$
is finite, there is an integer $c>0$ such that $S^{n-c}\cdot T\subset S^{n}.$
It follows that $|S^{n}\cap \phi ^{-1}(B)|\sim |T|\cdot |\overline{S}%
^{n}\cap \overline{\phi }^{-1}(B)|$ where $\overline{S}$ is the image of $S$
in $\Gamma /T.$ Moreover the $\overline{S}^{n}$'s is a family of nicely
growing subsets of $\Gamma /T$ according to Corollary \ref{Panpan}.
Therefore Theorem \ref{main} follows from Theorem \ref{eq}. 
\endproof%

\proof[Proof of Corollary \ref{coro}]%
Apply Theorem \ref{main} to $B=U$ and $B=V$ with $\phi =id$ and take the
ratio.%
\endproof%

We finally state one last result which generalizes Theorem \ref{eq} in an
obvious way (we assume here $G$ simply connected) and whose proof we only
sketch because it is entirely analogous to the proof of Theorem \ref{eq} we
just gave and presents no additional difficulties. Recall from Section \ref
{nicely} Lemma \ref{Lemm1} that $M_{\infty }=\lim_{t\rightarrow +\infty
}\delta _{\frac{1}{t}}M$ is a connected subgroup of $N_{0}$. Let $c_{M}>0$ 
be
the constant from Lemma \ref{Lemm2}. We have:

\begin{theorem}
\label{unifdist}For any bounded Borel subset $B$ of $G$ with $vol_{G}(\partial B)=0$ we have the following weak convergence of measures on 
$N$,
\begin{equation*}
\lim_{T\rightarrow +\infty }\frac{1}{T^{d(\Gamma )-d(G)}}\sum_{\gamma \in
\phi ^{-1}(B)} {\bf{\Delta }}_{\delta _{\frac{1}{T}}(\gamma )}=c_{M}\cdot
vol_{G}(B)\cdot vol_{M_{\infty }}
\end{equation*}
where ${\bf{\Delta}}_{x}$ is the Dirac mass at $x$.
\end{theorem}

\proof[Proof sketch]%
Let $f\geq 0$ be a compactly supported function on $N$, let $h=\mathbf{1}_{B}
$ and set $F_{t}(x)=\sum_{\gamma \in \Gamma }f(\delta _{\frac{1}{T}}(x\gamma
))h(\phi (x\gamma ))$ for $x\in N.$ Let $V_{t}=\int_{M}f(\delta _{\frac{1}{T}%
}(m))dvol_{M}(m).$ By Lemma \ref{Lemm2}, $V_{t}/t^{\deg _{N}(\frak{m}%
)}\rightarrow c_{M}\cdot \int fdvol_{M_{\infty }}.$ When this is positive, the unique ergodicity of 
$M$ on $N/\Gamma $ yields as above the weak convergence $\left\langle
F_{t}^{B},\psi \right\rangle /V_{t}\rightarrow vol_{G}(B)\cdot
\int_{N/\Gamma }\psi d\nu $. The point-wise convergence is derived in a
similar way.%
\endproof%

\subsection{Local limit theorem for random walk averages}

In this last paragraph we prove Corollary \ref{coro2}. The proof makes use
of the results obtained by Varopoulos and Alexopoulos about convolution
powers of measures on nilpotent groups. We first explain some terminology
and refer the reader to \cite{Alex} for the details.

Let $\mu =\sum_{s\in S}\mu (s)\delta _{s}$ be a symmetric probability
measure on $G$, whose finite support $S$ generates a dense subgroup $\Gamma $
of $G$. Let $\phi :N\rightarrow G\simeq N/M$ be the map given by Malcev's
rigidity as before. We now define the sub-Laplacian associated to $\mu $ on $%
N$, by setting 
\begin{equation*}
L_{\mu }=\sum_{i=n_{1}+1}^{n_{1}+n_{2}}b_{i}X_{i}+\frac{1}{2}\sum_{1\leq
i,j\leq n_{1}}a_{ij}X_{i}X_{j}
\end{equation*}
where the $X_{i}$'s are $N$-left invariant vector fields on $N$ that form an
adapted basis for the Lie algebra $\frak{n}$ (i.e. $\frak{n}=\oplus _{p\geq
1}m_{p}\,$and $C^{q}(\frak{n})=\oplus _{p\geq q}m_{p}$ where $C^{q}(\frak{n}%
)=span\{X_{i},d_{i}\geq q\}$ and $n_{p}=\dim m_{p}$). We let $(\delta
_{t})_{t>0}$ be the corresponding one-parameter group of dilations. The
coefficients are defined by 
\begin{eqnarray*}
a_{ij} &=&\int x_{i}x_{j}d\mu (x)\text{~if~}1\leq i,j\leq n_{1} \\
b_{i} &=&\int x_{i}d\mu (x)\text{~if~}n_{1}<i\leq n_{1}+n_{2}
\end{eqnarray*}
By definition, the heat kernel $(p_{t})_{t>0}$ associated to $\mu $ on $G$
is the solution to the heat equation for $L_{\mu }$ that is 
\begin{equation}
\frac{\partial p_{t}}{\partial t}=L_{\mu }p_{t}  \label{heat}
\end{equation}
Since $\Gamma $ is dense in $G$, the matrix $(a_{ij})_{1\leq i,j\leq n_{1}}$
is non-degenerate and equation (\ref{heat}) has a unique solution $%
(p_{t})_{t>0}$ such that $p_{t}$ is smooth and positive everywhere with
total integral over $N$ equal to $1$. We will also need to consider the
``heat kernel at infinity'' $(p_{t}^{0})_{t>0}$ defined by the same
equation, except that the $X_{i}$'s are replaced by the corresponding
left-invariant vector fields for the graded group structure $N_{0}$
associated to the group of dilations $(\delta _{t})_{t>0}$ (see \ref
{dilations}) The $(p_{t}^{0})_{t>0}$ enjoy the scaling property 
\begin{equation}
t^{\frac{d(\Gamma )}{2}}p_{t}^{0}(\delta _{\sqrt{t}}(x))=p_{1}^{0}(x).
\label{scale}
\end{equation}

Let $|\cdot |$ be a quasi-norm on $N$. We have:

\begin{theorem}
(\cite{Alex}, Corollary 1.19) For every $\alpha ,$ $0<\alpha <\frac{1}{2},$
there is $c_{\alpha }>0$ such that $\forall n\in \Bbb{N}$, $\forall \gamma
\in \Gamma ,$%
\begin{equation}
n^{\frac{d(\Gamma )}{2}}\left| \mu ^{*n}(\gamma )-p_{n}(\gamma )\right| \leq 
\frac{c_{\alpha }}{n^{\alpha }}\exp \left( \frac{-|\gamma |^{2}}{c_{\alpha
}\cdot n}\right)   \label{comp1}
\end{equation}
\end{theorem}

\begin{theorem}
(\cite{Alex}, Theorem 6.6) $\exists c>0,$ $\forall n\in \Bbb{N},$ $\forall x$%
\begin{equation}
n^{\frac{d(\Gamma )}{2}}\left| p_{n}(x)-p_{n}^{0}(x)\right| \leq \frac{c}{%
\sqrt{n}}  \label{comp2}
\end{equation}
\end{theorem}

\begin{theorem}
(\cite{Alex}, Corollary 6.5) Let $d_{i}=(X_{i})$. $\exists c>0,$ $\forall
k\in \Bbb{N}$, $\forall i_{1},...,i_{k}\geq 0,$ $\forall t>0,$%
\begin{equation}
\left| X_{i_{1}}\cdot ...\cdot X_{i_{k}}p_{t}(x)\right| \leq \frac{c}{%
t^{(d(\Gamma )+d_{i_{1}}+...+d_{i_{k}})/2}}\exp \left( \frac{-|x|^{2}}{%
c\cdot t}\right)   \label{comp3}
\end{equation}
\end{theorem}

\proof[Proof of Corollary \ref{coro2}]%
We have $\mu ^{*n}(B)=\sum_{\gamma \in \phi ^{-1}(B)}\mu ^{*n}(\gamma ).$
Fix $\alpha \in (0,\frac{1}{2})$ and let $\varepsilon >0$ with $\varepsilon
\cdot (d(\Gamma )-d(G))<\alpha .$ Then combining (\ref{comp1}) and (\ref
{comp3}) we get $n^{\frac{d(G)}{2}}\sum_{|\gamma |>n^{\frac{1}{2}%
+\varepsilon }}\mu ^{*n}(\gamma )\ll n^{\frac{d(G)}{2}}\sum_{k>n^{\frac{1}{2}%
+\varepsilon }}k^{d(\Gamma )}\cdot \exp (\frac{-k^{2}}{cn})$ which tends to $%
0$ as $n$ tends to $+\infty .$ The same holds for $p_{n}(\gamma )$ or $%
p_{n}^{0}(\gamma )$ in place of $\mu ^{*n}(\gamma )$. So when computing $n^{%
\frac{d(G)}{2}}\mu ^{*n}(B)$, we may as well restrict to counting points $%
\gamma \in B$ with $|\gamma |\leq n^{\frac{1}{2}+\varepsilon }.$ According
to Theorem \ref{eq} $\#\{\gamma \in \Gamma ,|\gamma |\leq n,\gamma \in \phi
^{-1}(B)\}=O(n^{d(\Gamma )-d(G)})$, hence combining this estimate with (\ref
{comp1}) we get, 
\begin{equation*}
n^{\frac{d(G)}{2}}\mu ^{*n}(B)=\sum_{|\gamma |\leq n^{\frac{1}{2}%
+\varepsilon },\gamma \in \phi ^{-1}(B)}n^{\frac{d(G)}{2}}p_{n}(\gamma
)+o(1)=\sum_{\gamma \in \phi ^{-1}(B)}n^{\frac{d(G)}{2}}p_{n}^{0}(\gamma
)+o(1)
\end{equation*}
By the scaling property (\ref{scale}) $n^{\frac{d(G)}{2}}p_{n}^{0}(\gamma
)=n^{-\frac{d(\Gamma )-d(G)}{2}}p_{1}^{0}(\delta _{\frac{1}{\sqrt{n}}%
}(\gamma )).$ Then by Theorem \ref{unifdist}
\begin{equation*}
n^{\frac{d(G)}{2}}\mu ^{*n}(B)=\frac{1}{n^{\frac{d(\Gamma )-d(G)}{2}}}%
\sum_{\gamma \in \phi ^{-1}(B)}p_{1}^{0}(\delta _{\frac{1}{\sqrt{n}}}(\gamma
))+o(1)=c_{M}\cdot vol_{G}(B)\cdot \int_{M_{\infty
}}p_{1}^{0}dvol_{M_{\infty }}+o(1)
\end{equation*}

\endproof


\begin{thebibliography}{99}
\bibitem{Alex}  G. Alexopoulos, \textit{Random walks on discrete groups of
polynomial volume growth}, Annals of Proba. vol \textbf{30}, n. 2, (2002),
p. 723-801.

\bibitem{ArK}  V.I. Arnol'd and A.L. Krylov, \textit{Uniform distribution of
points on a sphere and certain ergodic properties of solutions of linear
ordinary differential equations in a complex domain}, Dokl. Akad. Nauk SSSR 
\textbf{148}, (1963) p. 9--12.

\bibitem{AuB}  L. Auslander and J. Brezin, \textit{Uniform distribution in
solvmanifolds}, Advances in Math. \textbf{7}, 111--144 (1971).

\bibitem{Bab}  M. Babillot, \textit{Points entiers et groupes discrets, de
l'analyse aux syst\`{e}mes dynamiques}, in Panoramas et Synth\`{e}ses SMF\
Monographs, \textbf{13}, ``Rigidit\'{e}, groupe fondamental et dynamique'',
1-119, (2002).

\bibitem{Bass}  H. Bass, \textit{The degree of polynomial growth of finitely
generated nilpotent groups}, Proc. London Math. Soc. (3) \textbf{25} (1972),
603--614.

\bibitem{Bel}  A. Bellaiche, \textit{The tangent space in sub-Riemannian
geometry}, in Sub-Riemannian Geometry, edited by A. Bellaiche and J-J.
Risler,1--78, Progr. Math., \textbf{144}, Birkhauser, (1996).

\bibitem{Bre}  E. Breuillard, \textit{Random walks on Lie groups}, preprint,
www.math.polytechnique/\symbol{126}breuilla/part0gb.pdf

\bibitem{Bre2}  E. Breuillard,\textit{\ Geometry of locally compact groups
with polynomial growth and shape of large balls}, preprint 2007.

\bibitem{Bre3}  E. Breuillard, \textit{Local limit theorems and
equidistribution of random walks on the Heisenberg group}, Geom. Funct. Anal
(GAFA), \textbf{15} (2005) no 1, 49p.

\bibitem{DRS}  W. Duke, N. Rudnick, P. Sarnak, \textit{Density of integer
points on affine homogeneous varieties}, Duke Math. J. \textbf{71} (1993),
no. 1, 143--179.

\bibitem{EMM}  A. Eskin, \textit{Counting problems and semisimple groups},
Proceedings of the International Congress of Mathematicians, Vol. II
(Berlin, 1998). Doc. Math. 1998, Extra Vol. II, 539--552.

\bibitem{EMS}  A. Eskin, S. Mozes, N. Shah, \textit{Unipotent flows and
counting lattice points on homogeneous varieties}, Ann. of Math. (2) \textbf{%
143} (1996), no. 2, 253--299.

\bibitem{Goo}  R. W. Goodman, \textit{Nilpotent Lie groups: structure and
applications to analysis}, LNM \textbf{562}, Springer-Verlag (1976).

\bibitem{Gro}  M. Gromov, \textit{Carnot-Carath\'{e}odory spaces seen from
within}, in Sub-Riemannian Geometry, edited by A. Bellaiche and J-J. Risler,
79-323, Birkauser (1996).

\bibitem{Guiv}  Y. Guivarc'h, \textit{Croissance polynomiale et p\'{e}riodes
des fonctions harmoniques}, Bull. Sc. Math. France \textbf{101}, (1973), p.
353-379.

\bibitem{Guiv2}  Y. Guivarc'h, \textit{Equir\'{e}partition dans les espaces
homog\`{e}nes,} Th\'{e}orie ergodique (Actes Journ\'{e}es Ergodiques,
Rennes, 1973/1974), pp. 131--142. Lecture Notes in Math., Vol. 532,
Springer, Berlin, (1976).

\bibitem{Kaz}  D.A. Kazhdan, \textit{Uniform distribution on a plane}, Trudy
Moskov. Mat. Ob. \textbf{14}, (1965) , p. 299-305.

\bibitem{Led}  F. Ledrappier,\textit{\ Ergodic Properties of some linear
actions}, Journal of Mathematical Sciences, Vol. \textbf{105}, No. 2, 2001

\bibitem{LeP}  E. Le Page, \textit{Th\'{e}or\`{e}mes quotients pour
certaines marches al\'{e}atoires}, Comptes Rendus Acad. Sc., \textbf{279}, s%
\'{e}rie A, n. 2, (1974).

\bibitem{Pan}  P. Pansu, \textit{Croissance des boules et des g\'{e}od\'{e}%
siques ferm\'{e}es dans les nilvari\'{e}t\'{e}s}, Ergodic Theory Dynam.
Systems \textbf{3} (1983), no. 3, 415--445.

\bibitem{Rag}  M. S. Raghunathan, \textit{Discrete subgroups of Lie groups},
Springer Verlag (1972).

\bibitem{Sta}  A. Starkov, \textit{Dynamical systems on homogeneous spaces},
Transl. of Math. Mono., \textbf{190}, AMS Providence, xvi+243 pp (2000).

\bibitem{Wey}  H. Weyl, \textit{\"{U}ber die gleichverteilung von Zahlen
mod. Eins}, Math. Ann. \textbf{77} (1916), 313-352.
\end{thebibliography}
\end{document}